\documentclass[11pt,reqno,a4paper]{amsart}

\usepackage[a4paper,margin=25mm]{geometry}
\usepackage{amsmath,amssymb,amsthm,mathtools}
\usepackage{enumitem,mathrsfs}
\usepackage{hyperref}

\usepackage{tikz}
\usetikzlibrary{arrows.meta,calc,patterns,decorations.pathreplacing}

\usepackage{tcolorbox}
\newtcolorbox{roadstep}[2]{
  colback=white,colframe=black!30,colbacktitle=black!5,coltitle=black,
  boxrule=0.4pt,arc=1mm,left=7pt,right=7pt,top=6pt,bottom=6pt,
  toptitle=4pt,bottomtitle=4pt,before skip=5pt,after skip=5pt,
  fonttitle=\normalsize\sffamily\bfseries,
  title={Step~#1: #2}}

\definecolor{GuideInk}{HTML}{1F2937}
\definecolor{GuideMuted}{HTML}{566170}
\definecolor{GuideBlue}{HTML}{1E568C}
\definecolor{GuideRed}{HTML}{AD4E32}
\definecolor{GuideGreen}{HTML}{237454}
\definecolor{GuideGreenDark}{HTML}{16533B}
\definecolor{GuideTint}{HTML}{E0F0E5}
\definecolor{GuidePale}{HTML}{F7F8FA}

\newtheorem{theorem}{Theorem}[section]

\newtheorem{lemma}[theorem]{Lemma}

\newtheorem{thm}[theorem]{Theorem}
\newtheorem{prop}[theorem]{Proposition}
\newtheorem{conjecture}[theorem]{Conjecture}
\newtheorem{claim}[theorem]{Claim}

\theoremstyle{remark}
\newtheorem{rmk}[theorem]{Remark}

\numberwithin{equation}{section}

\newcommand{\R}{\mathbb R}\newcommand{\Q}{\mathbb Q}

\newcommand{\supp}{\operatorname{supp}}

\newcommand{\ZZ}{\mathbb Z}
\newcommand{\CC}{\mathbb C}
\newcommand{\Tr}{\operatorname{Tr}}

\newcommand{\disc}{\operatorname{disc}}

\newcommand{\Ec}{\mathcal{E}}
\newcommand{\Cc}{\mathcal{C}}
\newcommand{\Ac}{\mathcal{A}}
\newcommand{\Bc}{\mathcal{B}}
\newcommand{\Lc}{\mathcal{L}}
\newcommand{\Oc}{\mathcal{O}}
\newcommand{\Sc}{\mathcal{S}}

\newcommand{\Rc}{\mathcal{R}}

\newcommand{\N}{\mathbb{N}}
\renewcommand{\P}{\mathbb{P}}
\renewcommand{\S}{\mathbb{S}}
\newcommand{\Z}{\mathbb{Z}}

\newcommand{\ebf}{\mathbf{e}}

\newcommand{\ind}{\mathbf{1}}

\newcommand{\norm}[1]{\left\lVert#1\right\rVert}

\newcommand{\dxi}{\,\mathrm{d}\xi}
\newcommand{\ds}{\,\mathrm{d}s}
\newcommand{\dt}{\,\mathrm{d}t}
\newcommand{\dx}{\,\mathrm{d}x}

\DeclareMathOperator{\GL}{GL}
\DeclareMathOperator{\Norm}{N}
\DeclareMathOperator{\Vol}{Vol}
\DeclareMathOperator{\Covol}{Covol}
\DeclareMathOperator{\STRUCTURE}{\mathbf{Structure}}
\DeclareMathOperator{\PRODUCT}{\mathbf{Product}}

\def\eqalign#1{\null\,\vcenter{\openup\jot\mathsurround\dimen12
  \ialign{\strut\hfil$\textstyle{##}$&$\textstyle{{}##}$\hfil
      \crcr#1\crcr}}\,}

\title[Power loss for the local Mizohata--Takeuchi estimate on smooth hypersurfaces]{Power loss for the local Mizohata--Takeuchi estimate on smooth hypersurfaces}

\author[Cairo]{Hannah Cairo}
\address[Hannah Cairo]{Department of Mathematics, William E. Kirwan Hall, University of Maryland, College Park, MD, U.S.}
\email{hmcairo@umd.edu}

\author[Fenves]{Inbo Gottlieb Fenves}
\address[Inbo Gottlieb Fenves]{University of Wisconsin-Madison, Madison, WI, U.S.}
\email{gottliebfenv@wisc.edu}

\author[Kova\v{c}]{Vjekoslav Kova\v{c}}
\address[Vjekoslav Kova\v{c}]{University of Zagreb Faculty of Science, Department of Mathematics, Bijeni\v{c}ka cesta 30, 10000 Zagreb, Croatia.}
\email{vjekovac@math.hr}

\author[Oliveira]{Itamar Oliveira}
\address[Itamar Oliveira]{School of Mathematics, The Watson Building, University of Birmingham, Edgbaston, Birmingham, B15 2TT, England.}
\email{i.oliveira@bham.ac.uk, oliveira.itamar.w@gmail.com}

\author[Zakharov]{Dmitrii Zakharov}
\address[Dmitrii Zakharov]{Massachusetts Institute of Technology, Cambridge, MA, U.S.}
\email{zakhdm@mit.edu}

\author[Zhang]{Ruixiang Zhang}
\address[Ruixiang Zhang]{Department of Mathematics, Evans hall, University of California, Berkeley, Berkeley CA, U.S. and  School of Mathematics, 1 Einstein Drive, Institute for Advanced Study, Princeton NJ, U.S.}
\email{ruixiang@berkeley.edu, rzhang@ias.edu}



\begin{document}

\begin{abstract}
We prove power-type lower bounds for the local Mizohata--Takeuchi constants for $C^{2,\beta}$ hypersurfaces with $0<\beta\leq1$ whose second fundamental form is nonzero at some point. This result, in particular, disproves the so-called local Mizohata--Takeuchi conjectures for the paraboloid and the sphere and, in turn, presents an obstruction to a general weighted approach to the theory of Fourier restriction of those classical surfaces.
\end{abstract}

\maketitle

    \setcounter{tocdepth}{1}
    \tableofcontents
    \setcounter{tocdepth}{2}

\section{Introduction}

\noindent
\textbf{Background.} In the late 1970s, Stein \cite{S} raised the question of whether weighted $L^{2}$ estimates for Fourier multiplier operators, notably the ball and Bochner--Riesz multipliers, could be controlled by Kakeya or Nikodym-type maximal operators (see \cite{ZhangBR} for a more detailed exposition about this question). In the context of Fourier restriction theory, this leads to what is commonly known as \textit{Stein's conjecture}, which states that weighted $L^{2}$ norms of Fourier extension operators are controlled by maximal X-ray averages of the weight along lines normal to the underlying hypersurface.

More precisely, let $\Sigma$ be a $C^2$ hypersurface in $\R^{d+1}$, and define the \textit{Fourier extension operator}
\begin{equation*}
    \Ec_\Sigma f(x,t) = \int_\Sigma f(\xi) e(\xi \cdot (x,t)) \,\mathrm{d}\sigma(\xi),
\end{equation*}
where $\mathrm{d}\sigma$ is the surface measure of $\Sigma$ and $e(s) = e^{2\pi i s}$. For a locally integrable weight $w\geq 0$ and $\nu\in\mathbb{S}^{d}$, define
\begin{equation*}
        Xw(\nu) = \sup_{v \in \nu^\perp}
        \int_\R w(v+t\nu)\,\mathrm{d}t.
    \end{equation*}


Stein's conjecture reads as follows:

\begin{conjecture}[Local Stein conjecture]
\label{conj:localStein}
    For every $\varepsilon > 0$ and $R \geq 1$, we have
    \begin{equation}\label{eq:localStein}
       \int_{B_{R}}|\mathcal{E}_{\Sigma}f(x,t)|^{2}w(x,t) \,\mathrm{d}x\,\mathrm{d}t
        \lesssim_{\Sigma,\varepsilon} R^\varepsilon
        \int_\Sigma |f(\xi)|^2 Xw(N(\xi)) \,\mathrm{d}\sigma(\xi),
    \end{equation}
for all $f \in L^2(\Sigma)$ and nonnegative weights $w$ supported in the ball $B_R$ of radius $R$ centred at the origin, where $N(\xi)$ denotes a unit normal to $\Sigma$ at $\xi$.
\end{conjecture}

An immediate consequence of \eqref{eq:localStein} would be the \textit{Mizohata--Takeuchi} conjecture:

\begin{conjecture}[Local Mizohata--Takeuchi conjecture]
\label{conj:localMT} For every $\varepsilon > 0$ and $R \geq 1$, we have
    \begin{equation}\label{eq:localMT}
       \int_{B_{R}}|\mathcal{E}_{\Sigma}f(x,t)|^{2}w(x,t) \,\mathrm{d}x\,\mathrm{d}t
        \lesssim_{\Sigma,\varepsilon} R^\varepsilon
        \norm{f}_{L^2(\Sigma)}^2
        \norm{Xw}_{L^\infty(\mathbb{S}^{d})},
    \end{equation}
for all $f \in L^2(\Sigma)$ and nonnegative weights $w$ supported in $B_R$.
\end{conjecture}

The Mizohata--Takeuchi conjecture \ref{conj:localMT} has its origins in the study of well-posedness for Schr\"odinger-type evolution equations in \cite{Mizohata} and \cite{Takeuchi}. In the 1990s, Barcel\'o, Ruiz and Vega \cite{BRV} studied a similar problem by establishing certain weighted estimates for the Helmholtz equation, and obtained applications to both Fourier restriction and the well-posedness of Schr\"odinger equations with time-dependent drift terms. The formulation \eqref{eq:localMT} stems from \cite{BRV}.

The global versions of Conjectures \ref{conj:localStein} and \ref{conj:localMT} state that
    \begin{equation}\label{eq:globalStein}
       \int_{\mathbb{R}^{d+1}}|\mathcal{E}_{\Sigma}f(x,t)|^{2}w(x,t) \,\mathrm{d}x\,\mathrm{d}t
        \lesssim_{\Sigma}
        \int_\Sigma |f(\xi)|^2 Xw(N(\xi)) \,\mathrm{d}\sigma(\xi),
    \end{equation}
and
    \begin{equation}\label{eq:globalMT}
       \int_{\mathbb{R}^{d+1}}|\mathcal{E}_{\Sigma}f(x,t)|^{2}w(x,t) \,\mathrm{d}x\,\mathrm{d}t
        \lesssim_{\Sigma}
        \norm{f}_{L^2(\Sigma)}^2
        \norm{Xw}_{L^\infty},
    \end{equation}
for all $f \in L^2(\Sigma)$ and weights $w\geq 0$. These global variants were first disproved by one of the present authors in \cite{Cairo}. See also \cite{BGNO} for a phase-space formulation of \eqref{eq:globalStein} and \eqref{eq:globalMT}, and \cite{BKNO} for a different set of counterexamples to \eqref{eq:globalStein} and \eqref{eq:globalMT}.

Recent positive results include the global estimate for $C^2$ uniformly convex planar curves with semialgebraic indicator weights of bounded complexity due to Hao, Wang and one of the present authors \cite{HWZ}, and the $R^{1/2+\varepsilon}$ local bound for the hyperbolic paraboloid in $\mathbb{R}^3$ due to Kim and Wang \cite{KW}.

A further motivation for Stein's conjecture is its potential to establish a direct connection between two of the central problems in harmonic analysis: the Kakeya maximal conjecture and the Fourier restriction conjecture\footnote{OpenAI has recently announced proofs of the Kakeya maximal conjecture and Fourier restriction estimates for the sphere and paraboloid in the interior of the conjectured range \cite{OpenAIKakeya,OpenAILocalSmoothing}. These developments strengthen the motivation for the present paper, which explains limitations of the proposed weighted approach to Kakeya and restriction for the sphere and the paraboloid.}. The latter is known to imply the former (see \cite{Mat15}), but combining Conjecture \ref{conj:localStein} with the conjectured $L^{d+1}$ bound for the Kakeya maximal operator would give
\begin{equation*}
    \norm{\Ec_\Sigma f}_{L^{2(d+1)/d}(B_R)}
    \lesssim_{\Sigma,\varepsilon} R^\varepsilon
    \norm{f}_{L^{2(d+1)/d}(\Sigma)},
\end{equation*}
which is the local endpoint formulation of the Fourier restriction conjecture for smooth compact hypersurfaces with nonvanishing Gaussian curvature; see \cite{BN}. Thus, Stein's conjecture would provide a mechanism for transferring geometric control of concentration along tubes into optimal Fourier extension estimates.

\vspace{.3cm}
\noindent
\textbf{Failure of Conjecture \ref{conj:localMT}.} Let the \textit{$R$-local MT constant} be defined by
\begin{equation*}
    \Cc_\Sigma(R) = \sup_{\substack{f \in L^2(\Sigma)}}\sup_{\supp w \subseteq B_R} \frac{\norm{\Ec_\Sigma f}_{L^2(w)}^2}{\norm{f}_{L^2(\Sigma)}^2 \norm{Xw}_{L^\infty}}.
\end{equation*}
The local Mizohata--Takeuchi conjecture \ref{conj:localMT} can then be restated as
\begin{equation*}
    \Cc_\Sigma(R) \lesssim_\varepsilon R^\varepsilon.
\end{equation*}

Guth \cite{Guth} noted that axiomatic wave-packet decomposition alone cannot prove Conjecture \ref{conj:localMT}; this was fleshed out by Carbery, Iliopoulou and Wang \cite{CIW}, where the bound $\Cc_\Sigma(R) \lesssim_\varepsilon R^{d/(d+2)+\varepsilon}$ was established for strictly convex $C^2$ hypersurfaces with nonvanishing Gaussian curvature, alongside its optimality using only wave-packet axioms. Two of the present authors \cite{CZ} subsequently showed that there exists a set of $C^2$ hypersurfaces $\Sigma$ for which the reverse inequality $\Cc_\Sigma(R) \gtrsim_\varepsilon R^{d/(d+2)-\varepsilon}$ holds (see also \cite{Fenves}). Thus, for a generic $C^2$ hypersurface, the extent to which the Mizohata--Takeuchi conjecture holds has been settled entirely, up to the endpoint. On the other hand, both \cite{CZ} and \cite{Fenves} deal with \textit{random} $C^k$ surfaces, and thus the possibility that Conjecture \ref{conj:localMT} may hold for $C^\infty$-surfaces (and in particular, the paraboloid and the sphere) remained open; see the recent discussion in \cite[Sec.~6]{Zhang}.

Our goal is to prove that Conjecture \ref{conj:localMT} fails for \textit{every} sufficiently smooth connected hypersurface $\Sigma$ that is not contained in a hyperplane. Recall that a hypersurface is of class $C^{2,\beta}$ for some $0<\beta\leq1$ if, near each point and after a rigid change of coordinates, it becomes the graph of a $C^2$ function $\phi$ whose Hessian is H\"{o}lder continuous with exponent $\beta$. Thus, $\|\textup{D}^2\phi(\xi)-\textup{D}^2\phi(\eta)\|\lesssim |\xi-\eta|^\beta$ holds for all $\xi,\eta$ in the same chart, where the implicit constant is allowed to depend on the chart.

The main result of this paper is the following:

\begin{thm}\label{thm:powerlossMT}
    Take $0<\beta\leq1$ and suppose that $\Sigma$ is a $C^{2,\beta}$ hypersurface whose second fundamental form is nonzero at some point. Then for all $R \gg_\Sigma 1$ we have
    \begin{equation*}
        \Cc_\Sigma(R) \gtrsim_\Sigma R^{\alpha_\beta},
    \end{equation*}
    where $0 < \alpha_\beta \ll 1$ depends only on the smoothness parameter $\beta$ (and is, in particular, independent of the dimension and the curvature of $\Sigma$).
\end{thm}

The proof brings together ideas from harmonic analysis and algebraic number theory. We first reduce Theorem \ref{thm:powerlossMT} to the construction of finite sets with large parabolic additive energy and controlled line occupancy. The next step is to use a tensor-power argument to amplify the gain from the counterexample in \cite{Cairo} to an exponential one in very high dimensions. The last challenge is to transfer this gain back to the real line without destroying the underlying additive relations or the required separation properties, which can be done with tools from algebraic number theory. The resulting constructions yield the desired power loss for the parabola, which is then transferred to general curved hypersurfaces.

\vspace{.3cm}
\noindent
\textbf{Outline of the paper.} This paper is organised as follows. In Section \ref{sec:outline}, we present an outline of the argument. In Section \ref{sec:transfer}, we reduce Theorem \ref{thm:powerlossMT} to a statement about additive energy of certain well-structured configurations on the parabola. Section \ref{sec:arithmetic} uses techniques from algebraic number theory to produce power-large energy of these configurations, which are then applied to complete the proof of Theorem \ref{thm:powerlossMT} in Section \ref{sec:main}. Appendix \ref{app:number-fields} records the standard number-theoretic facts which we use in Section \ref{sec:arithmetic}.

\vspace{.3cm}
\noindent
\textbf{Notation.} For quantities $X,Y$, we write $X \lesssim Y$ or $X = O(Y)$ if $X \le CY$ for some absolute constant $C$; $X \sim Y$ means $X \lesssim Y \lesssim X$. We write $X \ll Y$ if $X \le cY$ for $c$ which is sufficiently small. If the implicit constants depend on some auxiliary parameter $D$, we write for example $X \lesssim_D Y$. The notation $\lessapprox$ will be introduced in Section \ref{sec:arithmetic}.
We also write $\# A$ for the cardinality of a finite set $A$.

\vspace{.3cm}
\noindent
\textbf{Acknowledgements.} The authors acknowledge the support of their individual fellowships and grants.

H.C. is supported by NSF GRFP.

I.F. was supported by the NSF RTG grant No. DMS-2230900.

V.K. was supported in part by the Croatian Science Foundation under the project HRZZ-IP-2022-10-5116 (\emph{FANAP}) and in part by the European Union -- NextGenerationEU through the National Recovery and Resilience Plan 2021--2026, via an institutional grant of the University of Zagreb Faculty of Science, IK IA 1.1.3, \emph{Impact4Math}.

I.O. was supported by his EPSRC Postdoctoral Fellowship UKRI3285. Part of this work was conducted while he was in residence at the Institut Mittag-Leffler in Djursholm, Sweden, during the autumn semester of 2026.

D.Z. was supported by the Simons Dissertation Fellowship.

R.Z. is supported by NSF CAREER DMS-2143989 and Sloan Research Fellowship.
He would like to thank Alex Cohen for helpful discussions.

\section{Outline of the argument}\label{sec:outline}

The outline of the proof of Theorem \ref{thm:powerlossMT} is organised into three parts: the analytic reduction
in Steps~(1)--(2), the combinatorial construction in Steps~(3)--(4),
and its arithmetic realisation in Steps~(5)--(7). We present the
ideas in the order that motivates them.
\vspace{.3cm}

\begin{enumerate}[label=(\arabic*)]
    \item \underline{\textit{Reduction to the parabolic case.}} The non-vanishing curvature assumption on $\Sigma$ ensures that, for some $2$-plane $\Pi$, the curve $\Gamma = \Sigma \cap \Pi$ is well approximated by a parabola (Lemma~\ref{lemma:localparaboloid}). Our function $f$ and weight $w$ will essentially lie on this one-dimensional slice, extended in the remaining dimensions isotropically. This allows us to reduce power loss on $\Sigma$ to power loss on the parabola $\P^1$; the H\"older assumption ensures that the approximation $\Gamma \approx \P^1$ is close enough to maintain power loss.
    \vspace{.3cm}
    \item \underline{\textit{Reduction to parabolic energy bounds.}} We perform a further reduction to show that the Mizohata--Takeuchi conjecture on $\P^1$ implies certain energy bounds on ``well-separated'' sets $A \subseteq \R$, $B \subseteq \R^2$ (Proposition~\ref{prop:transfer}). More precisely, it is enough to make the quotient
    \begin{equation}\label{eq:energyratiooutline}
        \frac{E(A,B)}{L\,\#A\,\#B}
    \end{equation}
    grow as a positive power of $H$, where $H^{-1}$ is the separation scale, $E(A,B)$ is the parabolic additive energy, and $L$ is the ``maximal strip occupancy'' of the set $B$. The energy gives a lower bound for the weighted extension norm (Lemma~\ref{lemma:weightedL^2norm}), while strip occupancy gives an upper bound for the X-ray norm of $w$ (Lemma~\ref{lemma:Xraynorm}). Our goal is then reduced to finding such sets $A,B$.
    \vspace{.3cm}
    \item \underline{\textit{Finding a seed configuration.}} We begin by finding ``seed'' sets $A_0\subseteq\mathbb{R}$ and $B_0\subseteq\mathbb{R}^{2}$ whose energy ratio exceeds any prescribed constant (see~\eqref{eq:seedgain}); these seeds will be combinatorial in nature, based on the counterexamples from \cite{Cairo}. This will translate to a disproof of the global Mizohata--Takeuchi conjecture \eqref{eq:globalMT}, but will not by itself give power loss.
    \vspace{.3cm}
    \item \underline{\textit{Amplification by tensor products.}} Next, we consider the large product sets $A_0^k \subseteq \R^k$, $B_0^{k} \subseteq \R^{2k}$; we refer to this product configuration as $\PRODUCT_0$ \eqref{eq:product0}. Then one can see that $\PRODUCT_0$ witnesses power loss for the desired energy ratio, with the relevant gain becoming exponential in the number $k$ of tensor factors, but does so in $\R^k$ instead of $\R$. This forces us to deal with two problems:
    \begin{enumerate}[label=(\roman*)]
        \item The coordinate projection maps $\R^k \to \R$ are highly non-injective, and directly projecting $\PRODUCT_0$ will lose all of the gain.
        \item The natural way to deal with (i) would be to slightly rotate our configuration so that the coordinate maps do become injective on $\PRODUCT_0$. This introduces the more serious issue that projection will no longer be a ring homomorphism, and in particular will not preserve the energy equation; moreover, we will have very little control on the point separation after projection.
    \end{enumerate}
    \vspace{.3cm}
    \item \underline{\textit{Introducing an auxiliary lattice.}} To deal with (ii), we note that, if we could find $\PRODUCT_0$ inside the ring of integers of a totally real number field $K$, embedded in $\R^k$ via its Minkowski embedding, then both difficulties (i) and (ii) above would be resolved; we refer to such a high-dimensional lattice as $\STRUCTURE$. Indeed, in such a case all coordinate projections are injective ring homomorphisms (see the discussion following~\eqref{eq:ABdef}). However, this is unlikely to work directly, since $\PRODUCT_0$ is discrete (and quite sparse), so there is a priori no reason to expect $\PRODUCT_0 \cap \STRUCTURE$ to be large.
    \vspace{.3cm}
    \item \underline{\textit{Thickening and intersecting.}} We next perform a parabolic rescaling, enlarging $\PRODUCT_0$ to be at a scale much larger than $\STRUCTURE$, and then thickening the points of this rescaled copy to be rectangular regions whose sizes are large relative to $\STRUCTURE$ but small relative to the rescaling; call the resulting object $\PRODUCT$ \eqref{eq:product}. Assuming the scales are chosen correctly, all of the thickenings in $\PRODUCT$ remain well separated (Lemma~\ref{lemma:seedrectangles}), but now we can expect that $\STRUCTURE$ will intersect these large rectangular regions (see Figure \ref{fig:thickening-intersection}). To make this precise, we require quantitative dimension-independent control on intersections of $\STRUCTURE$ with axis-aligned boxes (Lemma \ref{lemma:structureintersection}).
    \vspace{.3cm}
    \item \underline{\textit{Projection back to the real line.}} Once all of this is arranged, we choose the parameters in the following order. Fix the root discriminant $D$, then choose the seed size $m\gg_D 1$ and the rescaling parameter $\rho\gg_{D,m} 1$, independently of $k$. We then let $k=[K:\Q]$ tend to infinity through the powers of $2$ supplied by Lemma~\ref{lemma:towerdata}. The intersection of $\PRODUCT$ with $\STRUCTURE$ gives us sets $\Ac\subseteq K$ and $\Bc\subseteq K^2$ so that the ratio \eqref{eq:energyratiooutline} has size at least $(c_Dm)^k$ inside of $K$ when lines are restricted to being non-vertical in $K \times K$, where $c_D$ is some quantity depending on $D$ but not $m,\rho,k$ (Lemma~\ref{lemma:energylowerbound}). Choosing $m$ sufficiently large relative to $D$ makes this quotient at least $2^k$ (see the proof of Proposition~\ref{prop:amplification}).

    We then project (using field homomorphisms) and parabolically rescale to obtain $A\subseteq\mathbb{R}$ and $B\subseteq\mathbb{R}^{2}$ without changing their cardinalities or parabolic energy (see~\eqref{eq:ABdef} and the discussion following it). The required separation and thin-strip occupancy follow from the arithmetic gap estimates (Lemmas~\ref{lemma:gaps} and~\ref{lemma:stripoccupancy}). The separation parameter is $H=C^k$, with $C>1$ independent of $k$, while the corresponding spatial radius is $R\approx_\Sigma H^{1+2/\beta}$ (Proposition~\ref{prop:transfer}). Thus the exponential gain becomes a positive power of $R$.
\end{enumerate}


Here are the rough sizes of the involved quantities, once $D$ is fixed:
\begin{equation*}
    m = D^{O(1)}, \quad \quad \rho = \exp(O_D(m)), \quad \quad L = \exp(O_D(k)), \quad \quad  C = \rho^{O(1)}.
\end{equation*}

\section{A transfer argument}\label{sec:transfer}

In this section, our goal is to reduce the local Mizohata--Takeuchi conjecture \ref{conj:localMT} on an arbitrary $C^{2,\beta}$ hypersurface to an energy estimate on the parabola $\P^1$; this will correspond to Steps~(1)--(2) of the outline.

Let $P(a) = (a,a^2)$, and define for any ring $\Rc$ and sets $A \subseteq \Rc$, $B \subseteq \Rc \times \Rc$ the \textit{parabolic additive energy}
\begin{equation*}
    E(A,B) = E_\Rc(A,B) = \#\left\{ (a,a',\zeta,\zeta') \in A \times A \times B \times B : P(a) + \zeta = P(a') + \zeta' \right\}.
\end{equation*}
We will typically suppress the dependence on the ring $\Rc$ from the notation.

\begin{prop}\label{prop:transfer}
    Suppose $H \gg 1$ and $L \ge 1$. Suppose $A \subseteq [-1,1]$ and $B \subseteq [-1,1]^2$ are finite nonempty sets so that $A$ is $H^{-1}$-separated, and so that
    \begin{equation*}
        \# \{(b,\tau) \in B : \tau - 2sb \in I\} \le L
    \end{equation*}
    for every $s \in \R$ and interval $I$ of length $H^{-1}$. Then if $0 < \beta \le 1$ and $\Sigma$ is a $C^{2,\beta}$ hypersurface whose second fundamental form is nonzero at some point, we have
    \begin{equation*}
        \Cc_\Sigma(\kappa H^{1+2/\beta}) \gtrsim_\Sigma \frac{E(A,B)}{L\,\#A\,\#B},
    \end{equation*}
    for some $\kappa \gg_\Sigma 1$.
\end{prop}

\begin{rmk}
    Notice that the strip occupancy of the set $B$ is measured with respect to strips with `vertical width' $H^{-1}$, i.e. ones intersecting the vertical lines in an interval of length $H^{-1}$. This is in contrast to the X-ray transform in Conjecture \ref{conj:localMT} which more naturally corresponds to strips of (actual) width $H^{-1}$. This distinction is due to the parabolic rescaling step in the proof below, and allows us to weaken the restrictions on the desired set $B$.
\end{rmk}

From now on, implicit constants will be permitted to depend on $\Sigma$. In the local coordinates of Lemma
\ref{lemma:localparaboloid}, write
$\mathrm{d}\sigma(\xi,\phi(\xi))=J(\xi)\,\mathrm{d}\xi$,
where the Jacobian $J$ satisfies $J(\xi)\sim 1$. Moreover, we identify a
function $f$ on the parameter domain with the
function $f/J$ on $\Sigma$, so that
\[
    \Ec_\Sigma f(x,t)
    =\int_{\R^d}f(\xi)
    e\bigl(x\cdot\xi+t\phi(\xi)\bigr)\,\mathrm{d}\xi,
\]
with $\|f\|_{L^2(\Sigma)}
\sim \|f\|_{L^2(\R^d)}$.
We begin by fixing some parameters
(the implicit constants below will be chosen later,
but they are not important at this stage):

\begin{equation}\label{implicitctes}
    h \sim H^{-1/\beta}, \quad \quad T \sim \frac{H}{h^2}, \quad \quad Y \sim \frac{H}{h}, \quad \quad \varepsilon \sim \frac{h}{H}.
\end{equation}
Note that we have $\varepsilon \ll h \ll 1 \ll Y \ll T$.

We begin with the following lemma, which follows from the fact that $\Sigma$ is $C^{2,\beta}$ with somewhere nonvanishing curvature.

\begin{lemma}\label{lemma:localparaboloid}
    Up to an element of $\GL_{d+1}(\R) \ltimes \R^{d+1}$, we may locally describe $\Sigma$ in the form
    \begin{equation*}
        \Sigma =_{\mathrm{loc}} \left\{ (\xi,\phi(\xi)) : \xi \in [-1,1]^d \right\},
    \end{equation*}
    where $\phi(\xi) = \sum_{j=1}^n \lambda_j\xi_j^2 + r(\xi)$, with $n \ge 1$, $\lambda_j = \pm 1$, $\lambda_1 = 1$, and $|r(\xi)| \lesssim |\xi|^{2+\beta}$.
\end{lemma}

Note that in the coordinates of Lemma \ref{lemma:localparaboloid}, if we let $\Pi = \R\ebf_1 + \R\ebf_{d+1}$, then $\Sigma \cap \Pi$ is approximated by the standard parabola. This change of coordinates will affect all estimates by a bounded multiplicative factor depending only on the surface $\Sigma$ (which we will absorb into our asymptotic notation).

Next, we set $\chi(y) = (1-|y|)_+$ and $\eta(z) = \hat{\chi}(z) = \left( \frac{\sin(\pi z)}{\pi z} \right)^2$. For each $a \in A$, we set $Q_a = ha\ebf_1 + [-\varepsilon,\varepsilon]^d$. We then set
\begin{equation*}
    f = \varepsilon^{-d/2}\sum_{a \in A} \ind_{Q_a}.
\end{equation*}
Since $A$ is $H^{-1}$-separated, the set $\{ha\ebf_1 : a \in A\}$ is $\varepsilon$-separated, so we immediately obtain
\begin{equation}\label{lemma:fL^2norm}
    \norm{f}_{L^2}^2 \sim \#A.
\end{equation}
Next, we choose the weights. Define
\begin{equation*}
    w(x,t) = w_1(x_1,t)w_2(x_2,\ldots,x_d),
\end{equation*}
where
\begin{equation*}
    w_1(x_1,t) = \chi(x_1/Y)\chi(t/T)\left| \sum_{(b,\tau) \in B} e\bigl( hbx_1 + h^2\tau t \bigr) \right|^2 \quad \mathrm{and} \quad w_2(x_2,\ldots,x_d) = \prod_{j=2}^d \chi(x_j/Y).
\end{equation*}
We see that $f$ is essentially an isotropic thickening of a 1-dimensional function corresponding to the slice $\Sigma \cap \Pi$, and similarly $w$ is an isotropic thickening of a 2-dimensional weight. Note that $w$ is continuous and supported in $B_{O(T)} = B_{O(H^{1+2/\beta})}$.

\begin{lemma}\label{lemma:weightedL^2norm}
    We have
    \begin{equation*}
        \norm{\Ec_\Sigma f}_{L^2(w)}^2 \gtrsim T \cdot E(A,B).
    \end{equation*}
\end{lemma}
\begin{proof} Developing $\|\mathcal{E}_{\Sigma}f\|_{L^{2}(w)}^{2}$ gives
\begin{equation*}
    \|\mathcal{E}_{\Sigma}f\|_{L^{2}(w)}^{2}=\varepsilon^{-d}\sum_{(a,a',\zeta,\zeta')\in A^{2}\times B^{2}}\int_{Q_{a}}\int_{Q_{a'}}K_{\zeta,\zeta'}(\xi,\xi')\,\mathrm{d}\xi\,\mathrm{d}\xi',
\end{equation*}
where for $\zeta,\zeta' \in B$ the kernel $K_{\zeta,\zeta'}(\xi,\xi')$ is given by
    \begin{equation*}
        K_{\zeta,\zeta'}(\xi,\xi') = \biggl[ YT \eta\bigl(Y( \xi_1 -\xi_1' + h(b-b'))\bigr)\eta\bigl( T( \phi(\xi) - \phi(\xi') + h^2(\tau-\tau') ) \bigr) \biggr]\biggl[ Y^{d-1} \prod_{j=2}^d \eta\bigl( Y(\xi_j - \xi_j') \bigr) \biggr].
    \end{equation*}
Since $K_{\zeta,\zeta'}(\xi,\xi')\geq 0$, restricting the sum above to the tuples $(a,a',\zeta,\zeta') \in A^2 \times B^2$ that contribute to the parabolic energy gives
    \begin{equation}\label{ineq1-22102026}
        \norm{\Ec_\Sigma f}_{L^2(w)}^2 \gtrsim \varepsilon^{-d}\sum_{\substack{(a,a',\zeta,\zeta') \in A^2 \times B^2 \\ P(a) + \zeta = P(a') + \zeta'}} \iint_{Q_a \times Q_{a'}} K_{\zeta,\zeta'}(\xi,\xi')\dxi\dxi'.
    \end{equation}
Fix $(a,a',\zeta,\zeta')$ in the summand of \eqref{ineq1-22102026}. Write $\xi = ha\ebf_1 + u$, $\xi' = ha'\ebf_1 + u'$ and $\zeta = (b, \tau)$, $\zeta' = (b',\tau')$ and notice that direct substitution and Lemma \ref{lemma:localparaboloid} give
    \begin{align*}
        |Y(\xi_1 - \xi_1' + h(b-b'))| &\lesssim Y\varepsilon, \\
        |T(\phi(\xi) - \phi(\xi') + h^2(\tau-\tau'))| &\lesssim T(h\varepsilon + \varepsilon^2 + h^{2+\beta}), \\
        |Y(\xi_j-\xi_j')| &\lesssim Y\varepsilon,
    \end{align*}
    and so by appropriately choosing the implicit constants in \eqref{implicitctes} we can make all of these values $\ll 1$, in which case we obtain $K_{\zeta,\zeta'} \gtrsim Y^dT$ on $Q_a \times Q_{a'}$, hence
    \begin{equation*}
        \norm{\Ec_\Sigma f}_{L^2(w)}^2 \gtrsim \varepsilon^{-d}E(A,B)Y^dT\varepsilon^{2d} \gtrsim TE(A,B). \qedhere
    \end{equation*}
\end{proof}

Finally, we must bound the X-ray transform of $w$.

\begin{lemma}\label{lemma:Xraynorm}
    For any affine line $\ell \subseteq \R^{d+1}$, we have
    \begin{equation*}
        \int_\ell w \lesssim LT\,\#B.
    \end{equation*}
\end{lemma}
\begin{proof}
    We begin with a one-dimensional claim: if $Z > 0$ and $(\omega_\zeta)_{\zeta \in B}$ are arbitrary frequencies so that every interval $I$ of length $Z^{-1}$ contains at most $L$ such frequencies, then
    \begin{equation}\label{eq:onedimlemma}
        \int_{\R} \chi(t/Z)\left| \sum_\zeta c_\zeta e(t\omega_\zeta) \right|^2\dt \lesssim ZL\,\#B
    \end{equation}
    for any coefficients $(c_\zeta)_{\zeta \in B}$ of absolute value 1. Indeed, this follows from expanding the square and applying the decay of $\eta$ alongside appropriate truncations.

    Now, let $\ell = \{y + s(\sin(\theta)\nu,\cos(\theta)) : s \in \R\}$, where $\nu \in \S^{d-1}$. We split into two cases.

    \textbf{Case 1}: $|\theta| \le \pi/4$. Making a change of variables $t = y_{d+1} + s\cos \theta$ yields the bound
    \begin{equation*}
        \int_\ell w \ds \lesssim \int_\R \chi(t/T) \left| \sum_{\zeta \in B} c_\zeta e(t\omega_\zeta) \right|^2\dt,
    \end{equation*}
    where
    \begin{equation*}
        c_\zeta = e(hb(y_1 - y_{d+1}\nu_1\tan(\theta))), \quad \quad \omega_\zeta = h^2 \left(\tau + \frac{\nu_1\tan(\theta)}{h}b \right).
    \end{equation*}
    Note that the implicit constant is uniform in our range of $\theta$. Then $\omega_\zeta$ lying in an interval of length $T^{-1}$ is equivalent to $\tau - 2sb$ lying in an interval of length $\sim H^{-1}$, for $s = -\frac{\nu_1\tan(\theta)}{2h}$, and hence we may apply the bound \eqref{eq:onedimlemma} to obtain
    \begin{equation*}
        \int_\R \chi(t/T) \left| \sum_{\zeta \in B} c_\zeta e(t\omega_\zeta) \right|^2\dt \lesssim TL\,\#B.
    \end{equation*}

    \textbf{Case 2}: $|\theta| \ge \pi/4$. Then some coordinate of $\nu$, say $\nu_j$, has size $\gtrsim 1$, and so we can make a change of variables using the trivial bound on the associated exponential sum, yielding
    \begin{equation*}
        \int_\ell w \ds \lesssim \int_\R \chi(x_j/Y)\,(\#B)^2\dx_j \lesssim Y\,(\#B)^2 \lesssim TL\,\#B,
    \end{equation*}
    where the final bound comes from the observation $\#B \lesssim LH$. 
\end{proof}

\begin{proof}[Proof of Proposition \ref{prop:transfer}]
    Choose $f \in L^2(\Sigma)$ and $w \colon B_{O(H^{1+2/\beta})} \to [0,\infty)$ as above. Then combining~\eqref{lemma:fL^2norm} with
Lemmas~\ref{lemma:weightedL^2norm} and~\ref{lemma:Xraynorm} yields
    \begin{equation*}
        \Cc_\Sigma(O(H^{1+2/\beta})) \ge \frac{\norm{\Ec_\Sigma f}_{L^2(w)}^2}{\norm{f}_{L^2}^2\norm{Xw}_{L^\infty}} \gtrsim \frac{E(A,B)}{L\,\#A\,\#B},
    \end{equation*}
    as desired.
\end{proof}

\section{Arithmetic amplification}\label{sec:arithmetic}

In this section, we provide a construction of power-large energy using techniques from algebraic number theory; this will comprise Steps (3)--(6) of the outline. This argument follows ideas of the disproofs of the unit-distance conjecture \cite{Alon} and the sum-product conjecture \cite{BSSZ}, as presented in Bloom \cite{Bloom}. Our goal is to prove the following proposition.

\begin{prop}\label{prop:amplification}
    There exists some absolute $C > 1$ so that the following holds. For all $k \in 2^\N$, there exist sets $A \subseteq [-1,1]$ and $B \subseteq [-1,1]^2$ satisfying the following conditions:
    \begin{enumerate}[label=(\roman*)]
        \item $A$ is $C^{-k}$-separated.
        \item For all $s \in \R$ and all intervals $I$ of length $C^{-k}$, we have
        \begin{equation*}
            \# \left\{ (b,\tau) \in B : \tau - 2sb \in I \right\} \le L_k,
        \end{equation*}
        for some $L_k \le \exp(O(k))$.
        \item We have
        \begin{equation*}
            \frac{E(A,B)}{L_k\,\#A\,\#B} \ge 2^k.
        \end{equation*}
    \end{enumerate}
\end{prop}

\begin{figure}[h]
\centering
\begin{tikzpicture}[
    font=\small,
    line cap=round,
    line join=round,
    every node/.style={inner sep=1.5pt},
    axis/.style={draw=black!70,line width=.45pt},
    guide/.style={draw=black!45,line width=.35pt},
    measure/.style={draw=black!75,line width=.4pt},
    >={Latex[length=1.3mm,width=.9mm]}
]

\begin{scope}[shift={(1.8cm,2cm)},x=1.8cm,y=1cm]
    \node at (0,1.42) {$A\subset[-1,1]$};

    \draw[axis,->] (-1.12,0)--(1.13,0);
    \draw[black,line width=.6pt] (-1,0)--(1,0);

    \foreach \x in {-1,1}{
        \draw[axis] (\x,-.065)--(\x,.065);
        \node[anchor=north] at (\x,-.13) {$\x$};
    }
    \node[anchor=west] at (1.17,0) {$a$};

    \foreach \x in {-.86,-.56,-.23,.30,.65,.88}{
        \fill[black] (\x,0) circle[radius=.95pt];
    }

    \draw[guide] (-.23,.12)--(-.23,.52);
    \draw[guide] (.30,.12)--(.30,.52);
    \draw[measure,<->] (-.23,.47)--(.30,.47);
    \node[anchor=south] at (.035,.60) {$\ge H^{-1}$};
\end{scope}

\begin{scope}[shift={(9cm,1.9cm)},x=1.9cm,y=1.9cm]
    \node at (0,1.22) {$B\subset[-1,1]^2$};

    \fill[black!10]
        (-1,-.7)--(1,.7)--(1,.825)--(-1,-.575)--cycle;
    \draw[guide] (-1,-.7)--(1,.7);
    \draw[guide] (-1,-.575)--(1,.825);

    \foreach \b/\t in {
        -.86/.72, -.72/.32, -.80/-.80,
        -.48/.82, -.42/-.68, -.24/.50,
        -.10/-.78, .22/.82, .30/-.50,
        .40/-.14, .82/-.76, .84/.88, .84/.26,
        -.68/-.4135, -.14/-.0355, .52/.4265
    }{
        \fill[black] (\b,\t) circle[radius=.95pt];
    }

    \draw[axis] (-1,-1) rectangle (1,1);
    \draw[axis,->] (-1,-1)--(1.12,-1);
    \draw[axis,->] (-1,-1)--(-1,1.12);

    \node[anchor=west] at (1.16,-1) {$b$};
    \node[anchor=south] at (-1,1.16) {$\tau$};

    \foreach \x in {-1,1}{
        \node[anchor=north] at (\x,-1.08) {$\x$};
    }
    \foreach \y in {-1,1}{
        \node[anchor=east] at (-1.16,\y) {$\y$};
    }

    \draw[guide] (-1.04,-.700)--(-1.19,-.700);
    \draw[guide] (-1.04,-.575)--(-1.19,-.575);
    \draw[measure] (-1.15,-.700)--(-1.15,-.575);
    \node[anchor=east] at (-1.24,-.6375) {$H^{-1}$};

    \node at (0,-1.40) {$\#(B\cap S_{s,I})\le L_k$};
\end{scope}

\end{tikzpicture}
\caption{Schematic geometry in Proposition~\ref{prop:amplification},
with $H=C^k$. The shaded strip
$S_{s,I}=\{(b,\tau):\tau-2sb\in I\}$ has vertical width
$|I|=H^{-1}$; the occupancy bound holds for every $s\in\mathbb R$
and every interval $I$ of this length.}
\label{fig:amplification-geometry}
\end{figure}
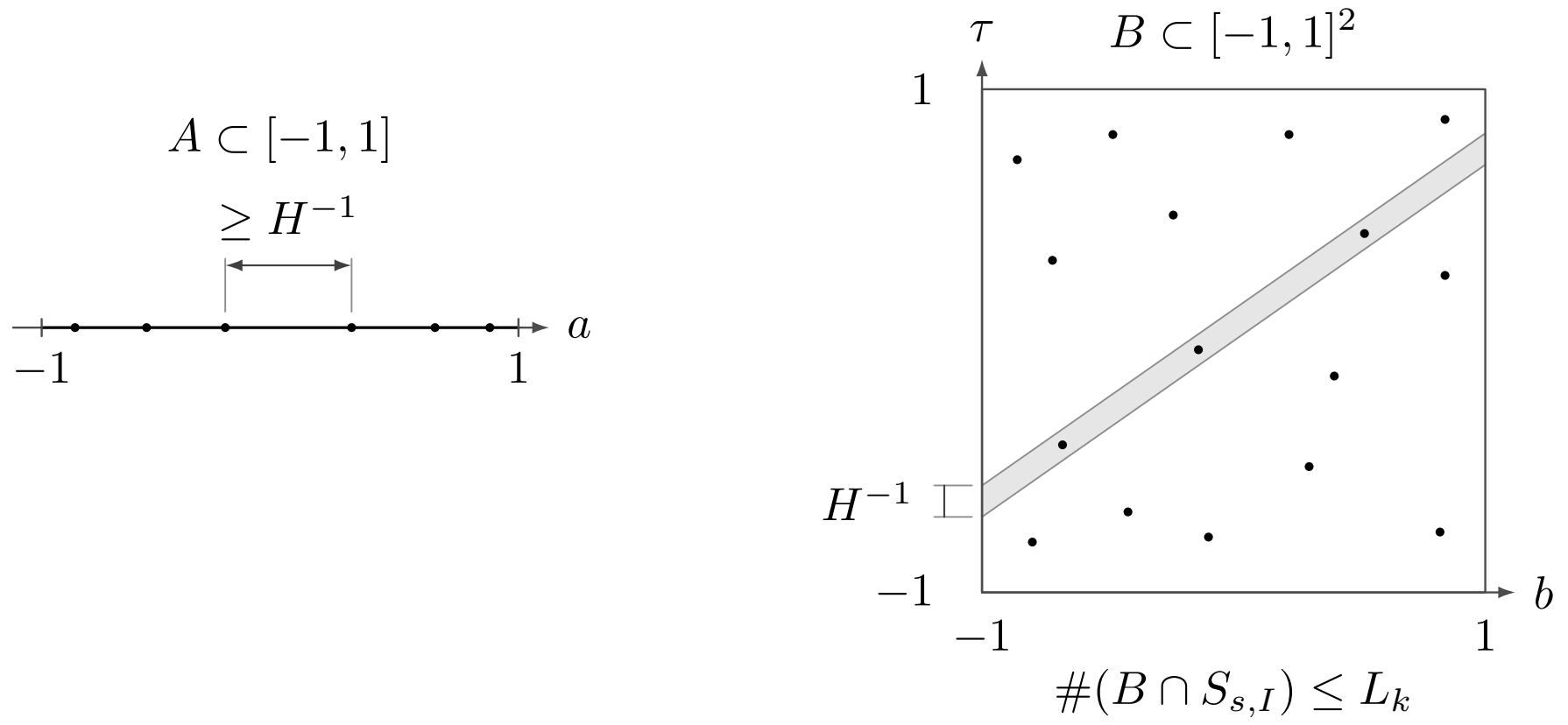

The geometry of conditions (i)--(ii) is illustrated in
Figure~\ref{fig:amplification-geometry}.

\subsection{\texorpdfstring{$\STRUCTURE$}{Structure}}

To start, we record the quantitative properties we need
of $\STRUCTURE$. Throughout this subsection, let $K/\Q$
be a totally real number field of degree $k$, with
embeddings $\sigma_1,\ldots,\sigma_k:K\hookrightarrow\R$ (see Appendix~\ref{app:number-fields} for other standard number-theoretic notation).
We will use the elementary norm bound from
Lemma~\ref{lemma:norminequalities}, proved in
Appendix~\ref{app:number-fields}. 

We require the following arithmetic input, which supplies us with an infinite family of totally real fields of bounded root discriminant. This construction is due originally to Martinet \cite{Martinet}.

\begin{lemma}\label{lemma:towerdata}
    There is some $D> 1$ so that, for all $j \in \N$, there exists a totally real field $K_j/\Q$ of degree $k = 2^{j+1}$ and root discriminant $D$.
\end{lemma}

Even though this will not affect in any way the results of this paper, it is perhaps helpful to think of $D$ as an explicit number. A concrete choice coming from Martinet's construction is
\[ D=4\sqrt{70035} =1058.565\ldots . \]
To simplify notation, we employ the following convention: once the root discriminant $D$ has been fixed as in Lemma \ref{lemma:towerdata}, write $X \lessapprox Y$ for $X \le b^kY$ for some $b > 1$ which is allowed to depend on $D$, but not on the parameters $m,\rho,k,C$ (some of which have not yet been introduced). Also write $X \approx Y$ if $X \lessapprox Y \lessapprox X$. We will always be able to take $m$ arbitrarily large in relation to $D$, and so this loss will always turn out to be tolerable.

\begin{lemma}\label{lemma:structureintersection}
    Let $K/\Q$ be a totally real number field of degree $k$ and root discriminant $D$, and let $\Lambda = \sigma(\Oc_K)$. Then
    \begin{equation*}
        \frac{\#(\Lambda \cap Q)}{\Vol(Q)/\Covol(\Lambda)} \approx 1
    \end{equation*}
    for any axis-aligned rectangle $Q = a + \prod_{i=1}^k [-M_i,M_i]$, with $M_i \ge C_0 D$ and a sufficiently large constant $C_0$ depending on $D$.
\end{lemma}
\begin{proof}
    The upper bound follows easily from a packing estimate, since distinct points in $\Lambda$ are $1$-separated in the $L^\infty$-metric, and $\Covol(\Lambda) = D^{k/2} \approx 1$. Indeed, for two distinct points $v,w\in\Lambda$ with $v=\sigma(\alpha),w=\sigma(\beta)$, we have $\alpha-\beta\in\Oc_{K}$. Hence its norm $\Norm_\Q^K(\alpha-\beta)\in\mathbb{Z}\backslash\{0\}$ and
    \[
    \prod_{i=1}^{k}|v_i - w_i|=|\Norm_\Q^K(\alpha-\beta)|\geq 1,
    \]
    which yields $\|v-w\|_{\infty}\geq 1$. This way, the cubes
\[
\widetilde{Q}_{v}:=v+\left[-\frac12,\frac12\right)^k,
\qquad v\in \Lambda\cap Q,
\]
are pairwise disjoint. Since $Q=a+\prod_{i=1}^k[-M_i,M_i]$, their union satisfies
\[
\bigcup_{v\in\Lambda\cap Q}\widetilde{Q}_{v}\ \subset\ a+\prod_{i=1}^k\left[-M_i-\frac12,M_i+\frac12\right].
\]
Therefore $\#(\Lambda\cap Q)
\lessapprox\Vol(Q)$. Since $\Covol(\Lambda)=D^{k/2}\approx 1$, this gives the desired upper bound. It thus suffices to show the lower bound.

To start, let $\Lambda^\vee$ be the dual lattice of $\Lambda$. Given any $\xi \in \Lambda^\vee$, we see that $\xi \cdot \Lambda \subseteq \Lambda^\vee$. Thus, if $\xi$ is nonzero, multiplication by $\xi$ sends a lattice of covolume $D^{k/2}$ to one of covolume at least $D^{-k/2}$ (see \eqref{covolumelambdavee}), and in particular the determinant of this transformation is at least $D^{-k}$; but this determinant is precisely\footnote{Here, we abuse notation slightly by identifying $K$ with $\sigma(K) \subseteq \R^k$.} $|\Norm_\Q^K(\xi)| = |\xi_1 \cdots \xi_k|$, and so
    \begin{equation}\label{eq:duallatticesep}
        |\xi_1\cdots \xi_k| \ge D^{-k} \gtrapprox 1
    \end{equation}
    for all $\xi \in \Lambda^\vee$ nonzero.

    Let $\psi \colon \R \to [0,1]$ be a fixed bump function supported on $[-1,1]$, with $\int \psi = 1$. Also assume we have Fourier decay $|\hat{\psi}(y)| \lesssim \min\{1,|y|^{-2}\}$. One can for example take $\psi \sim \ind_{[-1/2,1/2]} * \ind_{[-1/2, 1/2]}$. Choose a rectangle $Q = a + \prod_{i=1}^k [-M_i,M_i]$, let $\psi_i(x_i) = \psi(M_i^{-1}(x_i-a_i))$, and set $\psi_Q = \prod_{i=1}^k \psi_i$. By Poisson summation, we have
    \begin{equation*}
        \#(\Lambda \cap Q) \ge \sum_{v \in \Lambda} \psi_Q(v) = \frac{1}{D^{k/2}}\sum_{\xi \in \Lambda^\vee} \hat{\psi}_Q(\xi).
    \end{equation*}
    Note that $\hat{\psi}_Q(0) = M_1 \cdots M_k \hat{\psi}(0)^k = 2^{-k}\Vol(Q)$, and so it suffices to bound the remainder
    \begin{equation*}
        R := \sum_{\substack{\xi \in \Lambda^\vee \\ \xi \ne 0}} |\hat{\psi}_Q(\xi)|.
    \end{equation*}
    Fix integers $j_1,...,j_k \in \Z$, and consider the dyadic regions
    \begin{equation*}
        \Omega_{j_1,\ldots,j_k} = \{\xi \in \R^k : |M_i\xi_i| \in (2^{j_i},2^{j_i+1}]\}.
    \end{equation*}
    Since $|\hat{\psi}_i(z)| = M_i|\hat{\psi}(M_iz)|$, we have for $\xi \in \Omega_{j_1,\ldots,j_k}$ the estimate
    \begin{equation*}
        |\hat{\psi}_Q(\xi)| \lessapprox M_1 \cdots M_k \prod_{i=1}^k \min\{1,2^{-2j_i}\}.
    \end{equation*}
    Now, we estimate the intersections $\#(\Lambda^\vee \cap \Omega_{j_1,\ldots,j_k})$. For all $\xi',\xi'' \in \Lambda^\vee \cap \Omega_{j_1,\ldots,j_k}$, we have $|\xi_i' - \xi_i''| \lesssim M_i^{-1}2^{j_i}$. By pigeonholing in the first coordinate, we can moreover find some pair $\xi',\xi''$ for which
    \begin{equation}\label{eq:boundingdualintersect}
        |\xi_1' - \xi_1''| \lesssim \frac{M_1^{-1}2^{j_1}}{\#(\Lambda^\vee \cap \Omega_{j_1,\ldots,j_k})}.
    \end{equation}
    Applying \eqref{eq:duallatticesep} for $\xi = \xi'-\xi''$ and \eqref{eq:boundingdualintersect} and reorganising gives
    \begin{equation*}
        \#(\Lambda^\vee \cap \Omega_{j_1,\ldots,j_k}) \lessapprox M_1^{-1}\ldots M_k^{-1} 2^{j_1+\cdots +j_k} \approx \Vol(Q)^{-1} 2^{j_1+\ldots+j_k}.
    \end{equation*}
    Thus, denoting $j^+ = \max\{j,0\}$, we can estimate
    \begin{align*}
    R \lessapprox& M_1\ldots M_k\sum_{\overline{j}\in \ZZ^k} \#(\Lambda^\vee \cap \Omega_{\overline{j}})  2^{-2\sum_{i=1}^k j_i^+}
    \\\lessapprox& M_1\ldots M_k \Vol(Q)^{-1}\left(\sum_{j\in \ZZ} 2^{j - 2  j^+}\right)^k.
    \end{align*}
    The sum $\sum_{j\in \ZZ} 2^{j - 2  j^+} = \sum_{j <0} 2^j +\sum_{j\ge 0} 2^{-j}$ converges and so we obtain an overall estimate $R \lessapprox 1$. Recall that the main term contribution is $2^{-k} \Vol(Q) \approx \Vol(Q)$ and by the assumption $\Vol(Q) \ge (C_0 D)^k$. Taking $C_0$ larger than the constant in the notation $\lessapprox$, we conclude that the main term dominates and $\#(\Lambda \cap Q) \gtrapprox \Vol(Q) / \Covol(\Lambda)$.
\end{proof}

\begin{rmk}
    By carefully tracking through the constants in Lemma \ref{lemma:structureintersection}, one can show that the implicit constant $C_0$ and constants in the $\lessapprox$ notation are absolute (meaning independent of $D$), and in fact we have the bound
    \begin{equation*}
        2^{-(k+1)} \le \frac{\#(\Lambda \cap Q)}{\Vol(Q)/\Covol(\Lambda)} \le 2^{k+1}.
    \end{equation*}
\end{rmk}

\begin{rmk}
    Lemma \ref{lemma:structureintersection} can be viewed as a quantitative version of the well-known fact that the orbit of $\Lambda$ under the diagonal group $A = \{\mathrm{diag}(e^{t_1},\ldots,e^{t_k}) : t_1 + \ldots + t_k=0\}$ is compact in the space of lattices in $\R^k$. See \cite{SS} and the references therein for details on this perspective.
\end{rmk}

\subsection{\texorpdfstring{$\PRODUCT_0$ and $\PRODUCT$}
{Product 0 and Product}}

Next, we produce our base seed, based on \cite{Cairo}.

\begin{lemma}\label{lemma:seednoncollinearity}
For every integer $m\ge2$, there exist $m$ distinct integers
$a_1,\ldots,a_m\in\{1,\ldots,2^{3m}\}$ such that the points
\[
    \zeta_J=(b_J,\tau_J)
    =\sum_{j\in J}(a_j,a_j^2),
    \qquad J\subseteq\{1,\ldots,m\},
\]
are pairwise distinct and no three of them are collinear.
\end{lemma}

\begin{proof} For pairwise distinct subsets $J_1,J_2,J_3\subseteq\{1,\ldots,m\}$, consider
the polynomial
\[
    F_{J_1,J_2,J_3}(a_1,\ldots,a_m)
    :=
    \det\begin{pmatrix}
        b_{J_1}    & b_{J_2}    & b_{J_3} \\
        \tau_{J_1} & \tau_{J_2} & \tau_{J_3} \\
        1          & 1          & 1
    \end{pmatrix}.
\]
Writing $\varepsilon_{ir}=\ind_{\{i\in J_r\}}$,
multilinearity in the first two rows gives
\[
    F_{J_1,J_2,J_3}
    =\sum_{\substack{1\le i,j\le m\\i\ne j}}
       c_{ij}a_i^2a_j,
    \qquad
    c_{ij}
    =
    \det\begin{pmatrix}
        \varepsilon_{j1} & \varepsilon_{j2} & \varepsilon_{j3} \\
        \varepsilon_{i1} & \varepsilon_{i2} & \varepsilon_{i3} \\
        1 & 1 & 1
    \end{pmatrix}.
\]

We claim that $F_{J_1,J_2,J_3}$ is not identically zero.
Indeed, if that were the case, all $c_{ij}$ would vanish. Since the subsets $J_{1},J_2,J_3$ are pairwise distinct, at least one vector
$(\varepsilon_{i1},\varepsilon_{i2},\varepsilon_{i3})$ is neither $(0,0,0)$ nor $(1,1,1)$. Up to relabelling the subsets $J_i$, we may
assume without loss of generality that
$\varepsilon_{i2}=\varepsilon_{i3}\ne\varepsilon_{i1}$.
Computing $c_{ij}$ gives
\[
    0=c_{ij}
     =(\varepsilon_{i2}-\varepsilon_{i1})
      (\varepsilon_{j2}-\varepsilon_{j3})
    \qquad\text{for every }j.
\]
Since the first factor is $\pm1$, we obtain
$\varepsilon_{j2}=\varepsilon_{j3}$ for every $j$, which is a contradiction since $J_2\neq J_3$.

We now choose the $a_j$ so that all these determinants
are simultaneously nonzero. Fix a total order on the subsets of $\{1,\ldots,m\}$, and set
\[
    G(a_1,\ldots,a_m)
    :=\prod_{J_1<J_2<J_3}
       F_{J_1,J_2,J_3}(a_1,\ldots,a_m).
\]
Every factor is a nonzero polynomial of degree three, so
$G\neq0$ and $\deg G=3\binom{2^m}{3}<2^{3m}$. This way, $G$ cannot vanish at every point of
$\{1,\ldots,2^{3m}\}^m$.

Choose $(a_1,\ldots,a_m)$ in this grid with
$G(a_1,\ldots,a_m)\ne0$. By construction, no three points $\zeta_J$ generated from $\{a_{1},\ldots,a_{m}\}$ are
collinear. These points are also pairwise distinct:
if $\zeta_{J_1}=\zeta_{J_2}$ for two distinct subsets,
the determinant involving any third subset would be
zero. Finally, since the singleton points
$\zeta_{\{j\}}=(a_j,a_j^2)$ are distinct, the $a_j$
are distinct as well.
\end{proof}

\begin{rmk}
    One can in fact take $a_j = \ell^j$ in Lemma \ref{lemma:seednoncollinearity} for some $\ell \in \N_{\ge 2}$; the point of the above proof is that such configurations are \textit{generic} at lacunary scales.
\end{rmk}

Fix $a_1,\ldots,a_m$ as in
Lemma~\ref{lemma:seednoncollinearity}, and set
\[
    A_0=\{a_1,\ldots,a_m\},
    \qquad
    B_0=\{\zeta_J:J\subseteq\{1,\ldots,m\}\}.
\]
Then $\#A_0=m$, $\#B_0=2^m$, and no three points of
$B_0$ are collinear. Moreover, the seed coordinates
are integers of size $\exp(O(m))$, and the determinants
associated with triples of distinct points of $B_0$
have absolute value at least $1$.

For each $1\le j\le m$ and
$J\subseteq\{1,\dots,m\}$ with $j\notin J$, we have
$P(a_j)+\zeta_J=\zeta_{J\cup\{j\}}\in B_0$. Thus at least
$m2^{m-1}$ pairs $(a,\zeta)\in A_0\times B_0$ satisfy
$P(a)+\zeta\in B_0$. Let
$L_0=\sup_\ell \#(B_0\cap\ell)\le 2$, where the supremum is
over non-vertical affine lines. Define the \textit{representation function}

\begin{equation*}
    r(z):=\#\{(a,\zeta)\in A_{0}\times B_{0}: P(a)+\zeta=z\},
\end{equation*}
and observe that
\begin{equation}\label{rep-id-04112026}
    \begin{aligned}
    \displaystyle\sum_{z\in B_{0}}r(z)
    & = \displaystyle \sum_{z\in B_0}\sum_{a\in A_0}\sum_{\zeta\in B_0}\ind_{\{P(a)+\zeta=z\}}(a,\zeta) \\
    & = \displaystyle\sum_{a\in A_0}\#\{\zeta\in B_0:P(a)+\zeta\in B_0\} \\
    & \geq \displaystyle m2^{m-1}.
    \end{aligned}
\end{equation}
On the other hand, by Cauchy--Schwarz,
\begin{equation*}
    \left(\sum_{z\in B_0}r(z)\right)^{2}\leq \#B_0\,\left(\sum_{z\in B_0}r(z)^{2}\right) \leq \#B_0 E(A_0,B_0).
    \end{equation*}
Therefore
\begin{equation}\label{eq:seedgain}
        E(A_0,B_0)
        \ge \frac{(m2^{m-1})^2}{\#B_0}
         = \frac{m^2}{4}\#B_0\quad\Longrightarrow\quad
        \frac{E(A_0,B_0)}{L_0\,\#A_0\,\#B_0}\geq \frac{m}{8}.
\end{equation}
In particular, the seed energy ratio can be made arbitrarily
large by increasing $m$.

Taking the $k$-fold product, we define the
product configuration to be the ordered pair
\begin{equation}\label{eq:product0}
    \PRODUCT_0 := (A_0^k,B_0^k).
\end{equation}
Its two components are subsets of $\R^k$ and
$(\R^2)^k$, respectively. In the rest of this paper, we identify
$(\R^2)^k$ with $\R^k\times\R^k$ by grouping together the
horizontal and vertical coordinates:
\begin{equation}\label{identification-04112026}
     ((b_1,\tau_1),\ldots,(b_k,\tau_k))
    \longleftrightarrow
    ((b_1,\ldots,b_k),(\tau_1,\ldots,\tau_k)).
\end{equation}

Now, fix some $\rho\gg_{D,m}1$, in particular
requiring $\rho\ge\max\{2^{3m},2C_0 D,2\}$ (with $C_0$ from Lemma \ref{lemma:structureintersection}),
and define the sets
\begin{align*}
    I_j &= \rho a_j + [-C_0 D,C_0 D] \\
    R_J &= (\rho b_J + [-C_0 D,C_0 D]) \times (\rho^2\tau_J + [-\rho^{3/2},\rho^{3/2}]).
\end{align*}
Set
\[
    U:=\bigcup_{j=1}^m I_j,
    \qquad
    V:=\bigcup_{J\subseteq\{1,\ldots,m\}}R_J.
\]
The parabolically rescaled and thickened product
configuration is the ordered pair
\begin{equation}\label{eq:product}
    \PRODUCT := (U^k,V^k).
\end{equation}

\begin{lemma}\label{lemma:seedrectangles}
    For $\rho \gg_{D,m} 1$, we have $I_j \cap I_{j'} = \emptyset$ and $R_J \cap R_{J'} = \emptyset$ for each $j \ne j'$, $J \ne J'$. Moreover, any line in $\R^2$ intersects at most two rectangles in $V$.
\end{lemma}
\begin{proof}
    This all follows immediately from the choices of $A_0,B_0$ and from letting $\rho$ be sufficiently large. Note that $\rho$ will be of size $\exp(O(m))$ to guarantee the separation condition.
\end{proof}

\subsection{Intersecting
\texorpdfstring{$\PRODUCT$}{Product} with
\texorpdfstring{$\STRUCTURE$}{Structure}}

Now, we perform Step (6).

Intersecting the two components of $\PRODUCT$ with
$\Lambda$ and $\Lambda\times\Lambda$, respectively,
gives the pair
\begin{equation*}
    \Ac = \Lambda \cap U^k \quad \mathrm{and} \quad \Bc = (\Lambda \times \Lambda) \cap V^k. 
\end{equation*}
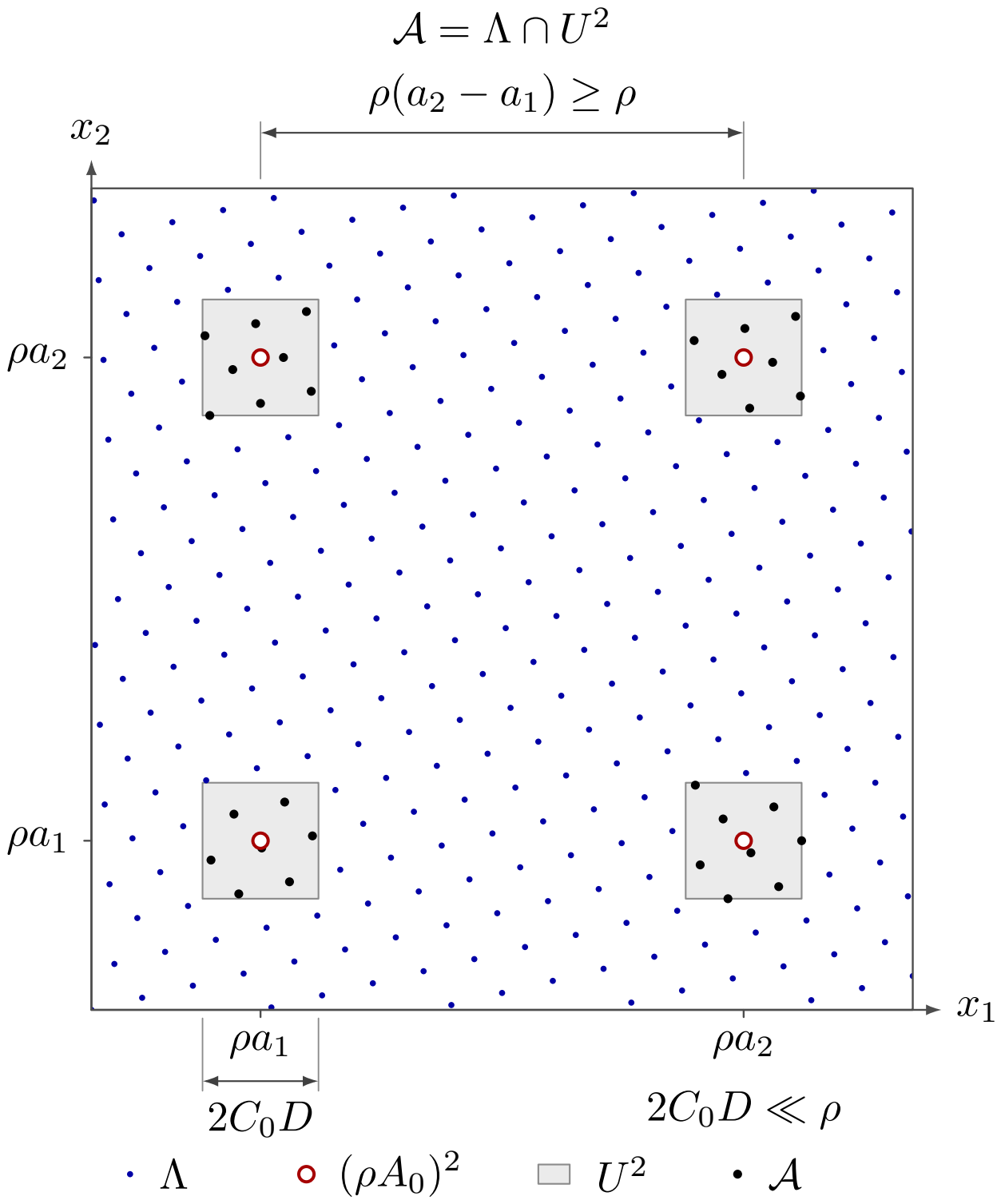
\begin{figure}[t]
\centering
\begin{tikzpicture}[
    x=1cm, y=1cm,
    font=\small,
    line cap=round, line join=round,
    every node/.style={inner sep=1.5pt},
    axis/.style={draw=black!70, line width=.45pt},
    guide/.style={draw=black!45, line width=.35pt},
    measure/.style={draw=black!75, line width=.4pt},
    box/.style={draw=black!45, fill=black!8, line width=.4pt},
    seed/.style={
        draw=red!65!black,
        fill=white,
        line width=.75pt
    },
    >={Latex[length=1.3mm,width=.9mm]}
]


\def\seedlow{1.4}
\def\seedhigh{5.4}

\def\boxhalf{.48}

\def\framesize{6.8}


\foreach \cx in {\seedlow,\seedhigh}{
    \foreach \cy in {\seedlow,\seedhigh}{
        \path[box]
            ({\cx-\boxhalf},{\cy-\boxhalf}) rectangle
            ({\cx+\boxhalf},{\cy+\boxhalf});
    }
}

%

\begin{scope}
    \clip (0,0) rectangle (\framesize,\framesize);

    \foreach \u in {-12,...,22}{
        \foreach \v in {-8,...,23}{

            \pgfmathsetmacro{\px}{.42*\u+.19*\v}
            \pgfmathsetmacro{\py}{.10*\u+.38*\v}

            \pgfmathtruncatemacro{\visiblepoint}{
                (\px>=0)&&(\px<=\framesize)&&
                (\py>=0)&&(\py<=\framesize)
            }

            \ifnum\visiblepoint=1\relax

                \pgfmathtruncatemacro{\selectedpoint}{
                    ((abs(\px-\seedlow)<=\boxhalf)||
                     (abs(\px-\seedhigh)<=\boxhalf))&&
                    ((abs(\py-\seedlow)<=\boxhalf)||
                     (abs(\py-\seedhigh)<=\boxhalf))
                }

                \ifnum\selectedpoint=1\relax
                    \fill[black]
                        (\px,\py) circle[radius=1.05pt];
                \else
                    \fill[blue!65!black]
                        (\px,\py) circle[radius=.72pt];
                \fi

            \fi
        }
    }
\end{scope}


\foreach \cx in {\seedlow,\seedhigh}{
    \foreach \cy in {\seedlow,\seedhigh}{
        \path[seed]
            (\cx,\cy) circle[radius=1.8pt];
    }
}


\draw[axis] (0,0) rectangle (\framesize,\framesize);

\draw[axis,->]
    (0,0)--({\framesize+.24},0);

\draw[axis,->]
    (0,0)--(0,{\framesize+.24});

\node[anchor=west]
    at ({\framesize+.30},0) {$x_1$};

\node[anchor=south]
    at (0,{\framesize+.30}) {$x_2$};

\foreach \c/\j in {\seedlow/1,\seedhigh/2}{
    \draw[axis] (\c,0)--(\c,-.06);
    \node[anchor=north]
        at (\c,-.12) {$\rho a_{\j}$};

    \draw[axis] (0,\c)--(-.06,\c);
    \node[anchor=east]
        at (-.12,\c) {$\rho a_{\j}$};
}


\draw[guide]
    (\seedlow,{\framesize+.08})
    --(\seedlow,{\framesize+.55});

\draw[guide]
    (\seedhigh,{\framesize+.08})
    --(\seedhigh,{\framesize+.55});

\draw[measure,<->]
    (\seedlow,{\framesize+.46})
    --(\seedhigh,{\framesize+.46});

\node[anchor=south]
    at ({(\seedlow+\seedhigh)/2},{\framesize+.54})
    {$\rho(a_2-a_1)\ge \rho$};

\node
    at ({\framesize/2},{\framesize+1.35})
    {$\mathcal A=\Lambda\cap U^2$};


\draw[guide]
    ({\seedlow-\boxhalf},-.07)
    --({\seedlow-\boxhalf},-.67);

\draw[guide]
    ({\seedlow+\boxhalf},-.07)
    --({\seedlow+\boxhalf},-.67);

\draw[measure,<->]
    ({\seedlow-\boxhalf},-.59)
    --({\seedlow+\boxhalf},-.59);

\node[anchor=north]
    at (\seedlow,-.70)
    {$2C_0D$};

\node
    at (\seedhigh,-.83)
    {$2C_0D\ll\rho$};


\fill[blue!65!black]
    (.32,-1.36) circle[radius=.72pt];
\node[anchor=west]
    at (.50,-1.36) {$\Lambda$};

\path[seed]
    (1.78,-1.36) circle[radius=1.8pt];
\node[anchor=west]
    at (1.98,-1.36) {$(\rho A_0)^2$};

\path[box]
    (3.70,-1.44) rectangle (3.96,-1.28);
\node[anchor=west]
    at (4.12,-1.36) {$U^2$};

\fill[black]
    (5.35,-1.36) circle[radius=1.05pt];
\node[anchor=west]
    at (5.54,-1.36) {$\mathcal A$};

\end{tikzpicture}

\caption{
A two-dimensional schematic of the parabolic rescaling and
thickening on the $\mathcal A$ side.
The hollow circles indicate the rescaled seed centres
$(\rho a_i,\rho a_j)$.
Each is thickened to a square of side length $2C_0 D$.
Since $2C_0D\ll\rho$, these rectangles remain well separated.
The oblique lattice $\Lambda$ is shown in blue, while the lattice
points lying in
$\mathcal A=\Lambda\cap U^2$ are shown in black.
}
\label{fig:thickening-intersection}
\end{figure}
For $u,v\in K$, let
\[
    \ell_{u,v}
    :=
    \{(b,\tau)\in K^2:\tau=ub+v\}
\]
be a non-vertical affine $K$-line, and write $\widetilde{\ell}_{u,v} :=  (\sigma\times\sigma)(\ell_{u,v})\subseteq\R^k\times\R^k$. We define
\[
    \Lc_K^{\mathrm{hor}}(\Bc)
    :=
    \sup_{u,v\in K}
    \#(\Bc\cap\widetilde{\ell}_{u,v}).
\]
Our goal is to prove the following energy estimate.

\begin{lemma}\label{lemma:energylowerbound}
    We have
    \begin{equation*}
        \frac{E_K(\Ac,\Bc)}{\Lc_K^{\mathrm{hor}}(\Bc)\,\#\Ac\,\#\Bc} \gtrapprox m^k.
    \end{equation*}
\end{lemma}
\begin{proof} We first estimate the sizes of $\Ac$ and $\Bc$.

\begin{claim}\label{claim1-04112026} The following bounds hold:
\begin{equation}\label{boundA}
    \#\Ac \gtrapprox m^k,
\end{equation}
\begin{equation}\label{boundB}
    \#\Bc \lessapprox (2^m\rho^{3/2})^k.
\end{equation}
\end{claim}
\begin{proof} The set $U^{k}$ is a disjoint union of $m^{k}$ axis-parallel boxes indexed by $\textnormal{\textbf{j}}=(j_{1},\ldots,j_{k})\in [m]^{k}$. Each such box is of the form $Q_{\textnormal{\textbf{j}}}=I_{j_1}\times\ldots\times I_{j_k}$, so $\Vol(Q_{\textnormal{\textbf{j}}})=(2C_0D)^{k}$. By Lemma \ref{lemma:structureintersection} we have
\begin{equation*}
    \#(\Lambda\cap Q_{\textnormal{\textbf{j}}})\gtrapprox \frac{\Vol(Q_{\textnormal{\textbf{j}}})}{\Covol(\Lambda)}=\frac{(2C_0D)^{k}}{D^{k/2}}\approx 1.
\end{equation*}
Since the boxes are pairwise disjoint,

\begin{equation*}
    \#\Ac=\sum_{\textnormal{\textbf{j}}\in [m]^{k}}\#(\Lambda\cap Q_{\textnormal{\textbf{j}}})\gtrapprox m^{k},
\end{equation*}
and \eqref{boundA} follows.

Similarly, $V^{k}$ is a disjoint union of $2^{mk}$ product boxes indexed by $\textnormal{\textbf{J}}=(J_{1},\ldots,J_{k})$, $J_i\subseteq\{1,\ldots,m\}$, and each has $k$ horizontal sides of length $2C_0D$ and $k$ vertical sides of length $2\rho^{3/2}$. The number of points in such a box is
\begin{equation*}
    \#\left(\Lambda\cap\prod_{i=1}^{k}[\rho b_{J_i}-C_0D,\rho b_{J_i}+C_0D]\right)\cdot\#\left(\Lambda\cap\prod_{i=1}^{k}[\rho^{2}\tau_{J_i}-\rho^{3/2},\rho^{2}\tau_{J_i}+\rho^{3/2}]\right).
\end{equation*}
Using the upper bound of Lemma \ref{lemma:structureintersection} and summing over all $2^{mk}$ boxes yields \eqref{boundB}.
\end{proof}

We next estimate the number of points of $\Bc$ that remain in $\Bc$
after translation by $P(\mathbf a)$, for each $\mathbf a\in\Ac$.

\begin{claim}\label{claim2-04112026}
For every $\mathbf a\in\Ac$, we have
\begin{equation*}
    \#\left\{
        \zeta\in\Bc:P(\mathbf a)+\zeta\in\Bc
    \right\}
    \gtrapprox\#\Bc,
\end{equation*}
where $P(\mathbf a)=(\mathbf a,\mathbf a^{\odot2})$
and $\mathbf a^{\odot2}$ denotes coordinatewise squaring.
\end{claim}

\begin{proof}
Fix $\mathbf a=\sigma(a)\in\Ac$, with $a\in\Oc_K$.
Since the intervals $I_j$ are disjoint, for each
$1\le i\le k$ there is a unique $j_i\in\{1,\ldots,m\}$
such that
\[
    x_i:=\sigma_i(a)=\rho a_{j_i}+u_i,
    \qquad |u_i|\le C_0D.
\]
For each choice of subsets
$J_i\subseteq\{1,\ldots,m\}\setminus\{j_i\}$, set
\[
    S_i:=R_{J_i}\cap
    \bigl(R_{J_i\cup\{j_i\}}-P(x_i)\bigr).
\]
Thus $(y,z)\in S_i$ precisely when
$(y,z)\in R_{J_i}$ and
$(y+x_i,z+x_i^2)\in R_{J_i\cup\{j_i\}}$.
In particular,
\begin{equation}\label{eq:overlap-membership}
    S_i\subseteq V,
    \qquad
    S_i+P(x_i)\subseteq V.
\end{equation}

We first check that these intersections are large.
The seed identities
$b_{J_i\cup\{j_i\}}=b_{J_i}+a_{j_i}$ and
$\tau_{J_i\cup\{j_i\}}=\tau_{J_i}+a_{j_i}^2$
show that the centres of the two rectangles defining
$S_i$ differ by
\[
    (-u_i,-2\rho a_{j_i}u_i-u_i^2).
\]
Since $a_{j_i}\le2^{3m}$, choosing $\rho\gg_{D,m}1$ gives $|2\rho a_{j_i}u_i+u_i^2|\leq \rho^{3/2}$. Both rectangles have width $2C_0D$ and height
$2\rho^{3/2}$, so $S_i$ has width at least $C_0D$
and height at least $\rho^{3/2}$.

We now count lattice points in the resulting product boxes.
There are $2^{m-1}$ choices of $J_i$ for each $i$, giving
$2^{(m-1)k}$ pairwise disjoint boxes
$S_1\times\cdots\times S_k$.
Indeed, $S_i\subseteq R_{J_i}$ and the rectangles $R_J$
are pairwise disjoint.
Under \eqref{identification-04112026}, each product is
$Q_b\times Q_\tau$, where $Q_b,Q_\tau\subseteq\R^k$
have side lengths at least $C_0D$ and $\rho^{3/2}$,
respectively.
Applying Lemma~\ref{lemma:structureintersection} to
both boxes, whose side lengths are at least $C_0 D$
for $\rho$ sufficiently large, gives
\[
\begin{aligned}
    \#((\Lambda\times\Lambda)\cap(Q_b\times Q_\tau))
    &=\#(\Lambda\cap Q_b)\,\#(\Lambda\cap Q_\tau)\\
    &\gtrapprox
    \left(\frac{(C_0D)^k}{D^{k/2}}\right)
    \left(\frac{\rho^{3k/2}}{D^{k/2}}\right) \\
    &\gtrapprox\rho^{3k/2}.
\end{aligned}
\]
Every point $\zeta=(\sigma(b),\sigma(\tau))$ counted
above belongs to $\Bc$, since its coordinate pairs
lie in $S_i\subseteq V$.
Moreover, \eqref{eq:overlap-membership} and
$x_i=\sigma_i(a)$ give
\[
    (\sigma_i(b+a),\sigma_i(\tau+a^2))
    =(\sigma_i(b)+x_i,\sigma_i(\tau)+x_i^2)
    \in V
    \qquad (1\le i\le k).
\]
Since $b+a,\tau+a^2\in\Oc_K$, this shows that
\[
    P(\mathbf a)+\zeta
    =(\sigma(b+a),\sigma(\tau+a^2))
    \in(\Lambda\times\Lambda)\cap V^k=\Bc.
\]

The product boxes are indexed by the tuples
$\mathbf J=(J_1,\ldots,J_k)$ with
$J_i\subseteq\{1,\ldots,m\}\setminus\{j_i\}$.
Since they are disjoint, summing their contributions gives
\[
\begin{aligned}
    \#\left\{\zeta\in\Bc:
        P(\mathbf a)+\zeta\in\Bc\right\}
    & \ge
    \sum_{\mathbf J}
    \#\left((\Lambda\times\Lambda)
        \cap(S_1\times\cdots\times S_k)\right) \\
    & \gtrapprox 2^{(m-1)k}\rho^{3k/2} \gtrapprox (2^m\rho^{3/2})^k \\
    & \gtrapprox\#\Bc,
\end{aligned}
\]
where the last inequality follows from \eqref{boundB}.
\end{proof}

We now prove the following energy lower bound.

\begin{claim}\label{claim3-04112026} $ E_K(\Ac,\Bc) \gtrapprox (\#\Ac)^2\,\#\Bc$.
\end{claim}
\begin{proof} Recall that
\begin{equation*}
    E_K(\Ac,\Bc)=\sum_{z\in\Lambda\times\Lambda}r(z)^{2},
\end{equation*}
where $r(z)$ is the representation function $r(z)=\#\{(a,\zeta)\in\Ac\times\Bc: P(\mathbf a)+\zeta=z\}$. Similarly to \eqref{rep-id-04112026},
\begin{equation*}
    \sum_{z\in\Bc}r(z)=\sum_{z\in\Bc}\sum_{a\in\Ac}\sum_{\zeta\in\Bc}\ind_{\{P(\mathbf a)+\zeta=z\}}(z)=\sum_{a\in\Ac}\#\{\zeta\in\Bc: P(\mathbf a)+\zeta\in \Bc\}\gtrapprox \#\Ac\,\#\Bc,
\end{equation*}
by Claim \ref{claim2-04112026}. By Cauchy--Schwarz, $(\sum_{z\in\Bc}r(z))^{2}\leq \#\Bc\,E_{K}(\Ac,\Bc)$, and the claim follows.
\end{proof}

It remains to bound the number of points of $\Bc$
on a non-vertical $K$-line.

\begin{claim}\label{claim4-04112026} $\Lc_K^{\mathrm{hor}}(\Bc) \lessapprox 1$.

\end{claim}
\begin{proof} Fix a non-vertical $K$-line
\[
    \ell=\{(b,\tau)\in K^2:\tau=ub+v\},
    \qquad u,v\in K,
\]
and consider its image under $\sigma\times\sigma$.
In the $i$th embedding, its points lie on the real line $y=\sigma_i(u)x+\sigma_i(v)$.

By Lemma~\ref{lemma:seedrectangles}, this line intersects
at most two of the rectangles $R_J$. Thus the points of
$\Bc$ in the image of $\ell$ lie in at most $2^k$
product boxes $R_{J_1}\times\cdots\times R_{J_k}$.
Under the identification~\eqref{identification-04112026},
each such product box has the form $Q_b\times Q_\tau$,
where $Q_b\subseteq\R^k$ has side lengths $2C_0D$.

Write $\widetilde{\ell}=(\sigma\times\sigma)(\ell)$.
Each $\sigma(b)\in\Lambda\cap Q_b$ determines at most
one point of $\Bc\cap\widetilde{\ell}$, since $\tau=ub+v$.
Thus Lemma~\ref{lemma:structureintersection} gives
\[
    \#(\Bc\cap\widetilde{\ell}\cap(Q_b\times Q_\tau))
    \le\#(\Lambda\cap Q_b)
    \lessapprox 1.
\]
Summing over the at most $2^k$ product boxes and taking
the supremum over all non-vertical $K$-lines, we obtain
\begin{equation*}
    \Lc_K^{\mathrm{hor}}(\Bc)
    =\sup_\ell\#(\Bc\cap\widetilde{\ell})
    \lessapprox 1. \qedhere
\end{equation*}
\end{proof}

Finally, by Claims \ref{claim1-04112026}, \ref{claim3-04112026} and \ref{claim4-04112026} we conclude that
\begin{equation*}
    \frac{E_K(\Ac,\Bc)}{\Lc_K^{\mathrm{hor}}(\Bc)\,\#\Ac\,\#\Bc}\gtrapprox \frac{(\#\Ac)^{2}\,\#\Bc}{\#\Ac\,\#\Bc}=\#\Ac\gtrapprox m^k,\end{equation*}
and the lemma follows.
\end{proof}

\subsection{Completing the proof}

To complete the proof of Proposition \ref{prop:amplification}, we need to project and parabolically rescale our sets $\Ac$ and $\Bc$, and then exploit the separation properties of $\Lambda$.

Choose some $S \gg_{\rho} 1$, and consider the sets
\begin{equation}\label{eq:ABdef}
    A=\left\{
        \frac{\sigma_1(a)}{S}:
        a\in\Oc_K,\ \sigma(a)\in\Ac
    \right\},
    \qquad
    B=\left\{
        \left(\frac{\sigma_1(b)}{S},
        \frac{\sigma_1(\tau)}{S^2}\right):
        b,\tau\in\Oc_K,\
        (\sigma(b),\sigma(\tau))\in\Bc
    \right\}.
\end{equation}
Note that for $S$ chosen sufficiently large, we have that $A \subseteq [-1,1]$ and $B \subseteq [-1,1]^2$. Since $\sigma_1$ is a field homomorphism, this also preserves both the cardinalities of the sets and their parabolic energy.

\begin{lemma}\label{lemma:gaps}
There exists $C\gg1$, depending on $D,m,\rho,S$
but not on $k$, such that every difference between
two elements of $A$, vertical difference between
two points of $B$, or determinant formed by three
points of $B$ is either zero or has absolute value
at least $100C^{-k}$.
\end{lemma}
\begin{proof}
By Lemma \ref{lemma:norminequalities}, if
$0\neq x\in\Oc_K$, then
\begin{equation}\label{ineq1-05112026}
    |\sigma_{1}(x)|\geq
    \frac{1}{\prod_{i=2}^{k}|\sigma_{i}(x)|}.
\end{equation}
Recall that the seed coordinates are bounded by $2^{O(m)}$
and $\rho$ was chosen sufficiently large in terms of
$D,m$. Thus $|\sigma_i(x)|\leq\rho^{100}$ whenever
$x$ is a difference of algebraic coordinates of $\Ac$
or $\Bc$, or a determinant formed from three algebraic
pairs corresponding to points of $\Bc$.

For distinct $a,a'\in\Oc_K$ with
$\sigma(a),\sigma(a')\in\Ac$, observe that
\begin{equation*}
    \left|\frac{\sigma_1(a)-\sigma_1(a')}{S}\right|
    =\frac{|\sigma_1(a-a')|}{S}
    \geq S^{-1}\rho^{-100(k-1)},
\end{equation*}
where the last inequality follows from
\eqref{ineq1-05112026} applied to $x=a-a'$.
We have analogous lower bounds for nonzero vertical
differences in $B$ and determinants of $B$-differences. Choosing $C=100S^3\rho^{100}$ completes the proof.
\end{proof}

Let us now introduce the notation
\begin{equation*}
    \Sc_k(B) = \sup_{s \in \R}\sup_{|I| = C^{-k}} \#\left\{ (b,\tau) \in B : \tau - 2sb \in I \right\}.
\end{equation*}
In words, $\Sc_k(B)$ is the maximum number of points of $B$ that a strip with vertical width $C^{-k}$ contains.

\begin{lemma}\label{lemma:stripoccupancy}
    We have $\Sc_k(B) \le \Lc_K^{\mathrm{hor}}(\Bc)$.
\end{lemma}

\begin{proof}
Fix $s\in\R$ and an interval $I$ of length $C^{-k}$,
and write
\[
    B_{s,I}:=\{(b,\tau)\in B:\tau-2sb\in I\}.
\]
By the determinant conclusion of Lemma \ref{lemma:gaps}, every three points of
$B_{s,I}$ are collinear. Moreover, any two distinct points
of $B_{s,I}$ have different horizontal coordinates. Thus, if $B_{s,I}$ contains at least two points,
they all lie on a single non-vertical real line $\ell$.

Let $(b_1,\tau_1),(b_2,\tau_2)\in\Oc_K^2$
be the algebraic pairs corresponding to two
distinct points of $B_{s,I}$. Since $b_1\neq b_2$,
we may define $u=(\tau_2 - \tau_1)(b_2 - b_1)^{-1}\in K$, $v=\tau_1 - ub_1\in K$. We can then lift $\ell$ to the $K$-line
\[
    \ell_K=\{(b,\tau)\in K^2:\tau=ub+v\}.
\]
By collinearity and injectivity of $\sigma_1$,
every point of $B_{s,I}$ lifts to a point of
$\Bc$ on $\ell_K$. Therefore
\[
    \#B_{s,I}\le\Lc_K^{\mathrm{hor}}(\Bc).
\]
The same inequality holds trivially if $B_{s,I}$
contains at most one point. Taking the supremum
over $s,I$ completes the proof.
\end{proof}

\begin{proof}[Proof of Proposition \ref{prop:amplification}]
    Fix $D$ as in Lemma \ref{lemma:towerdata}. Next, fix some $m \gg_D 1$ so large that the $m^k$ gain of Lemma \ref{lemma:energylowerbound} beats the $\gtrapprox 1$ by a definite factor; for example, choose $m$ so large that
    \begin{equation*}
        \frac{E_K(\Ac,\Bc)}{\Lc_K^{\mathrm{hor}}(\Bc)\,\#\Ac\,\#\Bc} \ge 2^k,
    \end{equation*}
    where $\Ac,\Bc$ are as in Lemma \ref{lemma:energylowerbound}, for every $k \in 2^\N$. Now, take $A,B$ as in \eqref{eq:ABdef}. By Lemma \ref{lemma:gaps}, we see that $A$ is $C^{-k}$-separated, completing (i). By Lemma \ref{lemma:stripoccupancy}, condition (ii) holds with $L_k = \Lc_K^{\mathrm{hor}}(\Bc)$. Since the slicing and parabolic rescaling preserved cardinalities and energy, Lemma \ref{lemma:energylowerbound} alongside our choice of $m$ implies (iii).
\end{proof}

\section{Proof of Theorem \ref{thm:powerlossMT}}\label{sec:main}

By combining Propositions \ref{prop:transfer} and \ref{prop:amplification}, we can complete the proof of Theorem \ref{thm:powerlossMT}.

\begin{proof}[Proof of Theorem \ref{thm:powerlossMT}]
    Let $C > 1$ be as in Proposition \ref{prop:amplification}. Let $R \gg_\Sigma 1$, and choose $k\in 2^\N$ such that
    \begin{equation*}
        \kappa C^{(1+2/\beta)k} \le R < \kappa C^{2(1+2/\beta)k},
    \end{equation*}
    where $\kappa$ is fixed by Proposition \ref{prop:transfer}. Now, set $H = C^k$, and $L = L_k$, where $C,L_k$ are defined as in Proposition \ref{prop:amplification}. Then the sets $A,B$ produced by Proposition \ref{prop:amplification} can be plugged into Proposition \ref{prop:transfer} to obtain
    \begin{equation*}
        \Cc_\Sigma(R) \geq \Cc_{\Sigma}(\kappa C^{(1+2/\beta)k})\gtrsim_{\Sigma} 2^k.
    \end{equation*}
    But note that $2^k > \kappa^{-\alpha_\beta}R^{\alpha_\beta}$, where
    \begin{equation}\label{powerloss}
        \alpha_\beta = \frac{\log 2}{2(1+2/\beta)\log C},
    \end{equation}
    which finishes the proof.
\end{proof}

For a $C^3$ hypersurface $\Sigma$ whose second fundamental form is nonzero at some point, explicit choices of parameters can be made in the previous proof to yield
\[ \Cc_\Sigma(R) \gtrsim_\Sigma R^{10^{-5}} \]
for all sufficiently large $R$. We made no attempt to optimise this exponent, in favour of expositional clarity; the exponent bound $10^{-5}$ is clearly improvable by a more careful choice of parameters. On the other hand, it is unlikely that a similar approach could match the best known upper bound $R^{d/(d+2)+\varepsilon}$ for every $\varepsilon>0$, proved by Carbery, Iliopoulou, and Wang \cite{CIW} for strictly convex compact patches in $\R^{d+1}$; indeed, there is no dependence on $d$ in \eqref{powerloss}.

We also note that the optimal power-loss exponent for a $C^2$ hypersurface $\Sigma \subseteq \R^{d+1}$, given by
\begin{equation*}
    \alpha_\Sigma := \limsup_{R \to \infty} \frac{\log \Cc_\Sigma(R)}{\log R},
\end{equation*}
depends on more than just the dimension $d$; indeed, for $C^2$ generic surfaces in $\R^3$ one has $\alpha_\Sigma = 1/2$, whereas \cite[Theorem 4.2]{Ortiz} implies that for the cone $\mathbb{C}\mathrm{one}^2 \subseteq \R^3$ (which is rank-1), one has $\alpha_{\mathbb{C}\mathrm{one}^2} \le 1/3$.


\appendix
\section{Number fields and lattices}\label{app:number-fields}

In this appendix we collect definitions and standard identities used in
Section~\ref{sec:arithmetic}; see, for instance,
\cite[Lectures~4, 6, 12 and~14]{Sutherland}.

\medskip
\noindent\textbf{Number fields and embeddings.}
A \emph{number field of degree $k$} is a field extension $K/\Q$
with $\dim_\Q K=k$. It has exactly $k$ distinct embeddings into
$\CC$, where an embedding is an injective field homomorphism
fixing $\Q$. The field is \emph{totally real} if all these
embeddings take values in $\R$. In this case, write them as
$\sigma_1,\ldots,\sigma_k\colon K\hookrightarrow\R$.
The \emph{Minkowski embedding} is
\[
    \sigma\colon K\longrightarrow\R^k,
    \qquad
    \sigma(x)=(\sigma_1(x),\ldots,\sigma_k(x)).
\]

\medskip
\noindent\textbf{Algebraic integers, norm and trace.}
An element of $K$ is an \emph{algebraic integer} if it is a root
of a monic polynomial with coefficients in $\Z$. These elements
form the \emph{ring of integers} $\Oc_K$. As an additive group,
$\Oc_K$ has a $\Z$-basis $\omega_1,\ldots,\omega_k$, called an
\emph{integral basis}, and $\Lambda=\sigma(\Oc_K)$ is a full
lattice in $\R^k$.
The \emph{field norm} and \emph{trace} of $x\in K$ are
\[
    \Norm_\Q^K(x)=\prod_{i=1}^k\sigma_i(x),
    \qquad
    \Tr_\Q^K(x)=\sum_{i=1}^k\sigma_i(x).
\]
Both belong to $\Q$, and both are integers when $x\in\Oc_K$.
In particular, the norm of a nonzero algebraic integer is a
nonzero integer.

\begin{lemma}\label{lemma:norminequalities}
For every nonzero $x\in\Oc_K$, we have
$|\Norm_\Q^K(x)|\ge1$. Moreover, if
$|\sigma_i(x)|\le T$ for all $1\le i\le k$,
with $T\ge1$, then
$|\sigma_i(x)|\ge T^{-k}$ for all $1\le i\le k$.
\end{lemma}

\begin{proof}
Since the norm of a nonzero algebraic integer is a
nonzero integer,
\[
    1\le|\Norm_\Q^K(x)|
    =\prod_{i=1}^k|\sigma_i(x)|.
\]
If $|\sigma_i(x)|\le T$ for every $i$, then,
for each $1\le j\le k$,
\[
    1\le|\sigma_j(x)|T^{k-1}.
\]
Hence $|\sigma_j(x)|\ge T^{-(k-1)}\ge T^{-k}$,
as claimed.
\end{proof}

\medskip
\noindent\textbf{Discriminant and root discriminant.}
The \emph{discriminant} of $K$ is
\[
    \disc K
    :=\det\bigl(\Tr_\Q^K(\omega_i\omega_j)\bigr)_{i,j=1}^k
    =\left(\det\bigl(\sigma_i(\omega_j)\bigr)_{i,j=1}^k\right)^2.
\]
It is independent of the integral basis and is a positive integer
because $K$ is totally real.
The \emph{root discriminant} is $D:=|\disc K|^{1/k}$.
The covolume of $\Lambda$, namely the volume of a fundamental
parallelepiped, satisfies
\[
    \Covol(\Lambda)=|\disc K|^{1/2}=D^{k/2}.
\]

\medskip
\noindent\textbf{The dual lattice.} For $\Lambda=\sigma(\Oc_K)$, the Euclidean dual lattice is
\[
    \Lambda^\vee
    :=\{\xi\in\R^k:\xi\cdot v\in\Z
        \text{ for every }v\in\Lambda\}.
\]
It satisfies $\Lambda^\vee\subseteq\sigma(K)$. Indeed, if
\[
    M=(\sigma_i(\omega_j))_{i,j=1}^k,
\]
then $\Lambda=M\Z^k$ and hence
$\Lambda^\vee=M^{-T}\Z^k$. Moreover,
\[
    G:=M^TM
    =\bigl(\Tr_\Q^K(\omega_i\omega_j)\bigr)_{i,j=1}^k
\]
has integer entries, so $G^{-1}$ has rational entries. Since $M^{-T}=M(M^TM)^{-1}=MG^{-1}$, we obtain $\Lambda^\vee\subseteq M\Q^k=\sigma(K)$.
Additionally, we have
\begin{equation}\label{covolumelambdavee}
     \Covol(\Lambda^\vee)=D^{-k/2}.
\end{equation}
Indeed, if we write $\Lambda=M\Z^k$ for some invertible real
matrix $M$, the defining condition for $\xi\in\Lambda^\vee$
is equivalent to $M^T\xi\in\Z^k$. Thus
$\Lambda^\vee=M^{-T}\Z^k$, and consequently
\[
    \Covol(\Lambda^\vee)
    =|\det M^{-T}|
    =\Covol(\Lambda)^{-1}
    =D^{-k/2}.
\]


\section*{Declaration of AI usage}
An early version of the proof of a special case of Theorem \ref{thm:powerlossMT} was generated by OpenAI's ChatGPT-5 Sol and ChatGPT-6 Astra from the previous works of the current authors. Namely, in a series of conversations held between 17 August and 10 September 2026 (\href{https://chatgpt.com/share/6a833050-3bb0-83eb-8004-4b48c00267af}{GPT-5 Sol chat}, \href{https://chatgpt.com/share/6a9d9c64-cb24-83eb-88cb-4742f1e75506}{GPT-6 Astra chat 1}, \href{https://chatgpt.com/share/6aa2a551-8fac-83ed-ba94-9be34bb9ea4b}{GPT-6 Astra chat 2}, \href{https://chatgpt.com/share/6aa31cdf-c3b8-83eb-84f3-2a8d2f3eddc3}{GPT-6 Astra chat 3}) those models were asked to gradually improve the lower bounds on the local Mizohata--Takeuchi constants from the paper \cite{BKNO}.\footnote{Interestingly, the authors also independently ran several experiments to investigate the chances of ChatGPT-6 Astra arriving at the proof completely autonomously, i.e., without strongly suggesting the context of \cite{BKNO}. Those chats did not arrive at a solution.} The prompts also directed the models to \cite{Cairo,CZ}. The original GPT construction relied on towers of number fields with prescribed prime splitting behaviour, provided by Hajir and Maire \cite{HM}. The authors then found a different, more direct proof presented here, which uses a different arithmetic construction in Section~\ref{sec:arithmetic} and generalises to $C^{2,\beta}$ hypersurfaces. ChatGPT-6 Astra was also used for copy-editing, for drafting the material in the short Appendix \ref{app:number-fields} on basic notions of algebraic number theory, and for proofreading of the final manuscript. The authors wrote the manuscript in their own words and assume full responsibility for its content.


\end{document}